\documentclass{llncs}
\title{On the Induction Operation for Shift Subspaces and Cellular Automata as Presentations of Dynamical Systems}
\author{Silvio Capobianco}
\institute{School of Computer Science, Reykjav\'{\i}k University
\\ Kringlan 1 - 103 Reykjav\'{\i}k, Iceland
\\ \email{silvio@ru.is} and \email{silvio.capobianco@gmail.com}
} %%end of \institute

\usepackage[english]{babel}
\usepackage{amsfonts}

\newcommand{\Acal}{\ensuremath{\mathcal{A}}}

\newcommand{\eg}{\textit{e.g.}}
\newcommand{\Forb}{\ensuremath{\mathcal{F}}}
\newcommand{\id}{\ensuremath{\mathrm{id}}}
\newcommand{\ie}{\textit{i.e.}}
\newcommand{\Neigh}{\ensuremath{\mathcal{N}}}

\newcommand{\restrict}[2]{\ensuremath{\left.{#1}\right|_{{#2}}}}

\newcommand{\supp}{\ensuremath{\mathrm{supp}\,}}

\newcommand{\vtheta}{\ensuremath{\vartheta}}
\newcommand{\Xset}{\ensuremath{\mathsf{X}}}

\newcommand{\Free}{\ensuremath{\mathbb{F}}}
\newcommand{\Nset}{\ensuremath{\mathbb{N}}}
\newcommand{\Zset}{\ensuremath{\mathbb{Z}}}

\begin{document}

\maketitle

\begin{abstract}
We consider continuous, translation-commuting transformations
of compact, translation-invariant families of mappings
from finitely generated groups into finite alphabets.
It is well-known that such transformations and spaces can be described
``locally'' via families of patterns and finitary functions;
such descriptions can be re-used on groups larger than the original,
usually defining non-isomorphic structures.
We show how some of the properties of the ``induced'' entities
can be deduced from those of the original ones, and vice versa;
then, we show how to ``simulate'' the smaller structure
into the larger one,
and obtain a characterization in terms of group actions
for the dynamical systems admitting of presentations via structures as such.
Special attention is given to the class of sofic shifts.
\end{abstract}

\section{Introduction} \label{sec_intro}
Cellular automata (briefly, CA)
are presentations of global dynamics in local terms:
the phase space is made of \emph{configurations}
on an underlying lattice structure,
and the transition function is induced by a pointwise evolution rule,
which changes the state at a node of the grid
by only considering finitely many \emph{neighbouring} nodes.
Originally, the only grids allowed were the hypercubic ones,
identified with the group $\Zset^d$ for some \emph{dimension} $d>0$,
and all the alphabets were finite,
though containing at least two elements;
this shall be referred to as the \emph{classical case}
in the rest of the present paper.

Modern CA theory borrows concepts and tools
from group theory, symbolic dynamics, and topology
(cf.~\cite{csms99,ff03,lm95}).
The lattice structure is provided by a \emph{Cayley graph}
of a finitely generated group:
the ``frames'' of this class generalize those of the classical case,
allowing more complicated grid geometries.
Such broadening, however, preserves
the requirement for finite neighbourhoods,
so that definition of global evolution laws in local terms
is still allowed.
Moreover, the phase space can be a \emph{subshift},
\ie, it can leave out some configurations,
but contains all of the \emph{translates} of each of its elements,
as well as the \emph{limits} of sequences it contains.
In general, however, the problem whether a given configuration
belongs to a given subshift, is only \emph{co-r.e.};
while this can be questionable when seeing CA as computation devices,
we cannot help but remark how the richer framework
simplifies dealing with \emph{simulations} between CA.

In this paper, which is an extended version
of a work submitted to the LATA 2008 conference~\cite{c08lata},
we deal with two problems.
The first one, is to understand when a dynamical system
can be described by a cellular automaton;
the second one, is to study the phenomena which happen
when a description for a subshift or a CA on a given group,
is employed in the context provided by a \emph{larger} group,
in the sense that the old one is a \emph{subgroup} of the new one.
At the time, we were not aware of the paper
by Ceccherini-Silberstein and Coornaert~\cite{csc08},
which also deals with induction of CA on larger groups,
and also considers a class of configuration spaces
which is broader than the classical one.
However, their work is focused on a broader class of \emph{alphabets}
without considering anything more general than the \emph{full shift}
that is made of all possible configurations
from a given group to a given alphabet;
on the other hand, our own work is aimed
towards the study of the most general subshifts,
provided the alphabet remains finite.

For the first problem, a solution is found employing \emph{group theory}:
a dynamical system admits of a CA presentation,
if and only if there exists a group action on it with special properties.
It is also observed, on one hand how the new class of CA
is strictly broader than the classical one;
and on the other hand, how some key properties of classical CA
are shared by the newer objects.

About the second one, a lemma about mutual inclusion
between images of shift subspaces via global CA functions,
showing that it is preserved \emph{either way}
when switching between the smaller group and the larger one.
This shall ensure that the operation of induction,
performed by ``recycling'' the description of the old object
(be it subshift or CA)
in the new \emph{context} given by the larger group,
is not only well defined,
but also independent on the specific description:
in other words, the induced object
only depends on the \emph{inducing} object.
We then show how several properties are transferred
from the old objects to the new ones,
some even either way as well;
this is of interest, because the new spaces and dynamics
is usually \emph{richer} than the old one.

A simulation of the original automaton into the induced one
is then explicitly constructed;
this extends to the case of arbitrary, finitely generated groups
the usual embedding of $d$-dimensional cellular automata
into $(d+k)$-dimensional ones.
This result puts another brick in the wall
of CA presentation of dynamical systems:
the class is not shrunk when the alphabet or the group are enlarged,
even up to bijections (for alphabets) and isomorphisms (for groups).
As a consequence of this fact,
the free group on two generators contains enough ``structure''
to present any CA dynamics on any free group.
Some remarks about \emph{sofic shifts}
are also made throughout the discussion.

The rest of the paper is organized as follows.
Section~\ref{sec_defs} provides a background.
Section~\ref{sec_cadyn} provides a characterization
for modern CA dynamics.
Section~\ref{sec_iss} provides the lemma of mutual inclusion
which ensures that the induction operation is well defined.
Section~\ref{sec_ica} studies induced CA
and how to embed the original CA into the induced one,
together with several considerations
for some special classes of subshifts.
Conclusions and acknowledgements follow.
%% end of Introduction

\section{Background} \label{sec_defs}
A \textbf{dynamical system} (briefly, d.s.) is a pair $(X,F)$
where the \textbf{phase space} $X$ is compact and metrizable
and the \textbf{evolution function} $F:X\to X$ is continuous.
If $Y\subseteq X$ is closed (equivalently, compact) and $F(Y)\subseteq Y$,
then $(Y,F)$ is a \textbf{subsystem} of $(X,F)$.
A \textbf{morphism} from a d.s. $(X,F)$ to a d.s. $(X',F')$
is a continuous $\vtheta:X\to X'$
such that $\vtheta\circ F=F'\circ\vtheta$;
an \textbf{embedding} is an injective morphism,
a \textbf{conjugacy} a bijective morphism.

Let $G$ be a group.
We write $H\leq G$ if $H$ is a subgroup of $G$.
If $H\leq G$ and $x\rho y$ iff $x^{-1}y\in H$,
then $\rho$ is an equivalence relation over $G$,
whose classes are called the \textbf{left cosets} of $H$,
one of them being $H$ itself.
If $J$ is a \textbf{set of representatives} of the left cosets of $H$
(one representative per coset)
then $(j,h)\mapsto jh$ is a bijection between $J\times H$ and $G$.

A \textbf{(right) action} of $G$ over a set $X$
is a collection $\phi=\{\phi_g\}_{g\in G}$
of transformations of $X$ (\ie, $\phi_g:X\to X$ for every $g\in G$)
such that
\begin{math}
\phi_{gh}=\phi_h\circ\phi_g
\end{math}
for all $g,h\in G$,
and $\phi_{1_G}=\id_X$, the identity function of $X$. 
Observe that the $\phi_g$'s are invertible,
with $(\phi_g)^{-1}=\phi_{(g^{-1})}$.
When $\phi$ is clear from the context, $\phi_g(x)$ can be written $x^g$.
Properties of functions (\eg, continuity) are extended to actions
by saying that $\phi$ has property $P$ iff each $\phi_g$ has property $P$.

If $G$ is a group and $S\subseteq G$,
the \textbf{subgroup generated by $S$}
is the set $\left<S\right>$ of all $g\in G$ such that
\begin{equation} \label{eq_gen}
g=s_1s_2\cdots s_n
\end{equation}
for some $n\geq 0$, with $s_i\in S$ or $s_i^{-1}\in S$ for all $i$.
$S$ is a \textbf{set of generators} for $G$ if $\left<S\right>=G$;
a group is \textbf{finitely generated} (briefly, f.g.)
if it has a finite set of generators (briefly, f.s.o.g.).
The \textbf{length} of $g\in G$ with respect to $S$
is the least $n\geq 0$ such that (\ref{eq_gen}) holds,
and is indicated by $\|g\|_S$.
The \textbf{distance} of $g$ and $h$ w.r.t. $S$
is the length $d^G_S(g,h)$ of $g^{-1}h$;
the \textbf{disk} of center $g$ and radius $R$ w.r.t. $S$
is $D^G_{R,S}(g)=\{h\in G\mid d^G_S(g,h)\leq R\}$.
In all such writings, $G$ and/or $S$ will be omitted
if irrelevant or clear from the context;
$g$, if equal to $1_G$.

An \textbf{alphabet} is a finite set with two or more elements;
all alphabets are given the discrete topology.
A \textbf{configuration} is a map $c\in A^G$
where $A$ is an alphabet and $G$ is a f.g. group.
Observe that the product topology on $A^G$
is induced by any of the distances $d_S$
defined by putting $d_S(c_1,c_2)=2^{-r}$,
$r$ being the minimum length w.r.t. $S$
of a $g\in G$ s.t. $c_1(g)\neq c_2(g)$.
Moreover, $\lim_{n\to\infty}c_n=c$ in the product topology
iff $\forall g\in G$ $\exists n_g\in\Nset$
such that $c_n(g)=c(g)$ for every $n>n_g$.

The \textbf{natural action} $\sigma^G$ of $G$ over $A^G$ is defined as
\begin{equation} \label{eq_sigmaG}
(\sigma^G_g(c))(h)=c(gh)\;\;\forall c\in A^G\;\;\forall g,h\in G\;;
\end{equation}
the superscript $G$ may be omitted if irrelevant or clear from the context.
Observe that $\sigma^G$ is continuous.
A closed subset $X$ of $A^G$ that is invariant by $\sigma^G$
is called a \textbf{shift subspace}, or briefly \textbf{subshift};
the case $X=A^G$ is called the \textbf{full shift}.
We use the notation $X\leq A^G$ to say that $X$ is a subshift of $A^G$.
The restriction of $\sigma^G$ to $X$
is again called the natural action of $G$ over $X$
and indicated by $\sigma^G$.
From now on, unless differently stated,
we will write $c^g$ for $\sigma^G_g(c)$.

Let $E\subseteq G$, $|E|<\infty$.
A \textbf{pattern} on $A$ with \textbf{support} $E$ is a map $p:E\to A$;
we write $E=\supp p$.
A pattern $p$ \textbf{occurs} in a configuration $c$
if there exists $g\in G$ such that $(c^g)_{|\supp p}=p$;
$p$ is \textbf{forbidden} otherwise.
Given a set $\Forb$ of patterns,
the set of all the configurations $c\in A^G$
for which all the patterns in $\Forb$ are forbidden
is indicated as $\Xset^{A,G}_\Forb$;
$A$ and/or $G$ will be omitted if irrelevant or clear from the context.
It is well known~\cite{ff03,lm95}
that $X$ is a subshift iff
$X=\Xset^{A,G}_\Forb$ for some $\Forb$.
$X$ is a \textbf{shift of finite type} if $\Forb$ can be chosen finite;
the full shift $A^G=\Xset^{A,G}_\emptyset$ is a shift of finite type.
A pattern $p$ is forbidden for $X\subseteq A^G$
if it is forbidden for all $c\in X$,
\ie $(c^g)_{|\supp p}\neq p$ for all $c\in X$, $g\in G$;
if $X$ is a subshift, this is the same as
$c_{|\supp p}\neq p$ for all $c\in X$.

A map $F:A^G\to A^G$ is \textbf{uniformly locally definable}
(UL-definable)
if there exist $\Neigh\subseteq G$, $|\Neigh|<\infty$,
and $f:A^\Neigh\to A$ such that
\begin{equation} \label{eq_ul-def}
(F(c))(g)=f\left(\restrict{c^g}{\Neigh}\right)
\end{equation}
for all $c\in A^G$, $g\in G$;
in this case, we write $F=F^{A,G}_f$.
Observe that any UL-definable function $F$ is continuous
and commutes with the natural action of $G$ on $A^G$;
\textbf{Hedlund's theorem}~\cite{ff03,hed69} states that,
if $X\subseteq A^G$ is a subshift and $F:X\to A^G$ is continuous
and commutes with the natural action of $G$ over $X$,
then $F$ is the restriction to $X$ of a UL-definable function.
Moreover, remark that, if $X$ is a subshift and $F$ is UL-definable,
then $F(X)$ is a subshift too:
if $X$ is of finite type, we say that $F(X)$ is a \textbf{sofic shift}.

A \textbf{cellular automaton} (CA)
with alphabet $A$ and \textbf{tessellation group} $G$
is a triple $\left<X_\Acal,\Neigh,f\right>$
where the \textbf{support} $X\subseteq A^G$ is a subshift,
the \textbf{neighbourhood index} $\Neigh\subseteq G$ is finite,
and the \textbf{local evolution function} $f:A^\Neigh\to A$
satisfies $F^{A,G}_f(X)\subseteq X$;
the restriction $F_\Acal$ of $F^{A,G}_f$ to $X$
is the \textbf{global evolution function},
and $(X_\Acal,F_\Acal)$ is the \textbf{associate dynamical system}.
Observe that $(X,F_\Acal)$ is a subsystem of $(A^G,F^{A,G}_f)$;
when $X=A^G$ is the full shift we say the CA is \textbf{full}.
Also observe that, because of Hedlund's theorem,
the class of CA with support $X$
can be seen as a monoid w.r.t. function composition.
When speaking of bijectivity, finiteness of type, etc. \emph{of $\Acal$},
we simply ``confuse'' it with either $F_\Acal$ or $X_\Acal$.
We say that $\Acal$ is \textbf{reversible}
if there exists a CA $\Acal'$,
\emph{with same alphabet, tessellation group, and support as $\Acal$},
such that $F_{\Acal'}\circ F_\Acal$ and $F_\Acal\circ F_{\Acal'}$
both coincide with the identity function \emph{of $X$}.
Observe that every reversible CA is bijective \emph{on its support}.

A pattern $p$ is a \textbf{Garden of Eden} (briefly, GoE)
for a CA $\Acal=\left<X,\Neigh,f\right>$
if it is allowed for $X$ and forbidden for $F_\Acal(X)$.
Any CA having a GoE pattern is nonsurjective;
compactness of $X$ and continuity of $F_\Acal$
ensure that the vice versa holds as well~\cite{ff03,mm93}.
$\Acal$ is \textbf{preinjective} if $F_\Acal(c_1)\neq F_\Acal(c_2)$
for any two $c_1,c_2\in X$ such that
\begin{math}
\{g\in G\mid c_1(g)\neq c_2(g)\}
\end{math}
is finite and nonempty.
\textbf{Moore-Myhill's theorem}~\cite{mo62,my62}
states that every full CA with tessellation group $\Zset^d$
is surjective iff it is preinjective.
This result has been extended
to larger classes of full CA~\cite{csms99,mm93},
but fails if the tessellation group
has a free subgroup on two generators~\cite{csms99}
or the support is not the full shift~\cite{ff03}

%% end of Background

\section{Characterization of CA Dynamics via Group Actions}
\label{sec_cadyn}

We have said in the intoduction that cellular automata
are presenytations of dynamical systems.
This remains only a nice, but vacuous concept
until we specify what it \emph{means}, for a CA,
to be a presentation:
intuitively, it should mean that the CA
``describes well'' the dynamics of the sistem.
How well, is stated in
\begin{definition} \label{def_CA-pres}
Let $(X,F)$ be a d.s., $\Acal$ a CA.
We say that $\Acal$ is a \emph{presentation} of $(X,F)$
if the latter and $(X,F_\Acal)$ are conjugate.
We call $CA(A,G)$ the class of d.s.
having a presentation as CA
with alphabet $A$ and tessellation group $G$.
We call $FCA(A,G)$ the subclass of $CA(A,G)$
made of d.s. having a presentation as CA
on the full shift $A^G$.
\end{definition}
One can wonder whether the introduction of CA on ``partial'' subshifts
is a factual extension of the concept.
Why should it not be possible to rewrite a system
using every possible configuration,
instead of only a selected package?
Why should we lose the feature of \emph{computability}
and step into the realm of \emph{recursive enumerability},
which is \textit{a priori} the only ensured thing
when leaving the full shift to accept arbitrary shift spaces?
Once we sell our soul to the devil of uncomputability,
we cannot get it back.

There can be many reasons to accept this kind of Faustian pact.
The first one, is that the new model looks promising
with respect to \emph{embeddings}:
it seems convenient to keep calling ``cellular automaton''
a local model of a subsystem of the associate dynamics.
Indeed, as we will see in the next sections,
the larger model actually behaves very well in this respect.

Another possible---and perhaps more compelling---reason, however,
is that the \emph{cardinality} of the phase space could be ``wrong'',
in the sense that it may hamper a presentation as a full CA;
however, if the system still displays the ``correct'' features,
we may still want to get a presentation in local terms,
and keep on calling it ``cellular automaton''.

And this is the content of
\begin{proposition} \label{prop_FCA-neq-CA}
Let $G$ be a f.g. group with $|G|\geq 3$.
\begin{enumerate}
\item \label{item:prop_FCA-neq-CA-fin}
If $G$ is finite, then there exists $X\leq A^G$
such that $|X|$ is not a perfect power.
\item \label{item:prop_FCA-neq-CA-inf}
If $G$ is infinite, then $A^G$ has a countable subshift.
\end{enumerate}
\end{proposition}
\begin{proof}{}
Fix $a,b\in A$ with $a\neq b$.

If $G$ is finite, the subset $X$ of configurations such that
\begin{itemize}
\item $c(g)\in\{a,b\}$ for every $c\in G$, and
\item there exist $g_a,g_b\in G$ s.t. $c(g_a)=a$ and $c(g_b)=b$
\end{itemize}
is closed and translation invariant, and has
\begin{math}
2^{|G|}-2=2\cdot(2^{|G|-1}-1)
\end{math}
elements, which is not a perfect power since $|G|\geq 3$.

If $G$ is infinite, then it is countable.
Let $X$ be the set of configurations such that
\begin{itemize}
\item $c(g)\in\{a,b\}$ for every $c\in G$, and
\item $c(g)=b$ for at most one $g\in G$.
\end{itemize}
Then $X$ is countable and translation invariant;
it is closed as well, because if $\lim_{n\to\infty}c_n=c$
and $c(g)$ is either outside $\{a,b\}$
or equals $b$ for two values $g_1,g_2\in G$,
then all the $c_n$'s for $n$ large enough
either take a value outside $\{a,b\}$
or take the value $b$ at $g_1$ and $g_2$,
and cannot belong to $X$.
\qed
\end{proof}
\begin{corollary} \label{cor_FCA-neq-CA}
If $|G|\geq 3$ then $CA(A,G)\neq FCA(B,H)$
for any alphabet $B$ and f.g. group $H$.
\end{corollary}
\begin{proof}{}
Because of Proposition~\ref{prop_FCA-neq-CA},
there exists $X\leq A^G$
which is not in bijection with a full shift.
Then no element of $CA(A,G)$ with support $X$
can be conjugate to an element of $FCA(B,H)$.

An immediate example of such CA is the identical transformation of $X$.
\qed
\end{proof}
\begin{example} \label{ex_CAnonF}
With the notations and conventions of Corollary~\ref{cor_FCA-neq-CA},
a less trivial example of a CA
whose associate d.s. admits of no presentation as a full CA,
can be constructed by fixing $\nu\in G\setminus\{1_G\}$,
putting $f_\nu(\nu\mapsto x)=x$,
and putting $\Acal_\nu=\left<X,\{\nu\},f_\nu\right>$,
where $X$ is as in the thesis of the corollary.
\qed
\end{example}
Regarding Example~\ref{ex_CAnonF},
it must be noted (cf.~\cite{ff03}),
that $F_{\Acal_\nu}$ is \emph{not}, in general,
the translation $\sigma_\nu$.
Actually, if $c\in A^G$, then
\begin{displaymath}
(F_{\Acal_\nu}(c))(g)
=f_\nu(\restrict{c^g}{\{\nu\}})
=c(g\nu)\;;
\end{displaymath}
to have this coincide with $c(\nu g)=(\sigma_\nu(c))(g)$
\emph{for every $g$ and $c$},
we must have $g\nu=\nu g$ for every $g\in G$,
that is, $\nu$ must belong to the \emph{center} of $G$.
This is a phenomenon already observed by Fiorenzi~\cite{ff03};
since it is useful to keep it in mind, we state it as
\begin{proposition} \label{prop_transl-ul}
Let $g\in G$.
Then $\sigma^G_g$, as a homeomorphism of $A^G$,
is UL-definable iff $g$ is central in $G$.
\end{proposition}
\begin{proof}{}
If $\sigma^G_g:A^G\to A^G$ is UL-definable,
then it commutes with $\sigma^G_h$ for every $h\in G$
because of Hedlund's theorem.
This implies, by evaluating the translates
$\sigma_{gh}(c)$ and $\sigma_{hg}(c)$ in $1_G$,
that $c(gh)=c(hg)$ for all $c\in A^G$, $h\in G$:
which is only possible if $gh=hg$ for all $h\in G$.
\qed
\end{proof}
We have thus given some reasons why
to deal with the more general concept of cellular automaton,
instead of sticking to the classical one.
Now that we know what a CA presentation \emph{is},
we must understand what a CA presentation \emph{requires}.

It turns out that the key feature of the dynamical systems
that admit of some presentation as CA,
also allow the tessellation group chosen for the CA
to act on their phase space \emph{like they were acting on a subshift}.
Here, ``like'' means that some key properties of $\sigma^G$
are shared 

We therefore state
\begin{definition} \label{def_discern}
Let $X$ be a set, $A$  an alphabet, $G$ a group,
$\phi$ an action of $G$ over $X$.
\emph{$X$ is discernible on $A$ by $\phi$}
if there exists a continuous function $\pi:X\to A$
such that, for any two distinct $x_1,x_2\in X$,
there exists $g\in G$ such that $\pi(\phi_g(x_1))\neq\pi(\phi_g(x_2))$.
\end{definition}
Observe, in Definition~\ref{def_discern}, the continuity requirement,
which demands that $\pi(x)=\pi(y)$
if $x,y\in X$ are ``near enough''.
\begin{example} \label{ex_discern}
Let $\Acal=\left<X,\Neigh,f\right>$ be a CA,
and let $(X,F_\Acal)$ be its associate d.s.
Then $\sigma_G$ commutes with $F_\Acal$.
Let $\pi(c)=c(1_G)$:
then, for any two $c_1,c_2\in X$,
$c_1(g)\neq c_2(g)$ is the same as
\begin{math}
\pi(\sigma^G_g(c_1))\neq\pi(\sigma^G_g(c_2)).
\end{math}
\qed
\end{example}
From Example~\ref{ex_discern} we know that
any CA dynamics admits of a ``discerning action''---which
happens to just be the natural action.
This gives us the hope that we have got a good clue
about the properties that characterize CA dynamics.

And this is confirmed by
\begin{theorem} \label{thm_char}
Let $A$ be an alphabet, $G$ a f.g. group, $(X,F)$ a d.s.
The following are equivalent:
\begin{enumerate}
\item \label{thm_char_p_present}
$(X,F)\in CA(A,G)$;
\item \label{thm_char_p_action}
there exists a continuous action $\phi$ of $G$ over $X$
such that $F$ commutes with $\phi$
and $X$ is discernible on $A$ by $\phi$.
\end{enumerate}
\end{theorem}
\begin{proof}
We start with supposing that $\Acal=\left<X_\Acal,\Neigh,f\right>$
is a presentation of $(X,F)$.
Let $\theta:X\to X_\Acal$ be a conjugacy from $(X,F)$ to $(X_\Acal,F_\Acal)$;
put
\begin{displaymath}
\phi_g=\theta^{-1}\circ\sigma^G_g\circ\theta
\end{displaymath}
for all $g\in G$, and
\begin{displaymath}
\pi(x)=(\theta(x))(1_G)\;.
\end{displaymath}
Remark that $\phi=\{\phi_g\}_{g\in G}$ is an action of $G$ over $X$
and that $(\theta(x))(g)=(\theta(x))^g(1_G)$ for all $x$ and $g$.
Continuity of $\phi$ and commutation with $F$ are straightforward to verify.
If $x_1\neq x_2$,
then $(\theta(x_1))(g)\neq(\theta(x_2))(g)$ for some $g\in G$,
thus
\begin{displaymath}
\pi(\phi_g(x_1))
=(\sigma^ G_g(\theta(x_1)))(1_G)
\neq(\sigma^ G_g(\theta(x_2)))(1_G)
=\pi(\phi_g(x_2))\;.
\end{displaymath}
For the reverse implication,
let $\pi$ as in Definition~\ref{def_discern}:
then $\tau:X\to A^G$ defined by
\begin{displaymath}
(\tau(x))(g)=\pi(\phi_g(x))
\end{displaymath}
is injective.
Moreover,
\begin{math}
(\tau(\phi_g(x))(h)=\pi(\phi_h(\phi_g(x)))
=\pi(\phi_{gh}(x))=(\tau(x))(gh)
\end{math}
for every $x\in X$, $g,h\in G$:
thus, $\tau\circ\phi_g=\sigma^G_g\circ\tau$ for all $g\in G$,
and $X'=\tau(X)$ is invariant under $\sigma^G$.

We now prove that $\tau$ is continuous.
Let $\lim_{n\in\Nset}x_n=x$ in $X$:
by continuity of $\pi$ and $\phi$,
$\lim_{n\in\Nset}(\tau(x_n))(g)=(\tau(x))(g)$ in $A$ for all $G$.
Since $A$ is discrete, for each $g\in G$ there exists $n_g$
such that $\pi(\phi_g(x_n))=\pi(\phi_g(x))$ for every $n>n_g$:
this is the definition of convergence of $\tau(x_n)$ to $\tau(x)$
in the product topology of $A^G$.

Since $X$ and $A^G$ are compact and Hausdorff,
$X'$ is closed in $A^G$ and a subshift,
while $\tau$ is a homeomorphism between $X$ and $X'$.
Define $F':X'\to X'$ by $F'=\tau\circ F\circ\tau^{-1}$:
then $(X',F')$ is a d.s.
and $\tau$ is a conjugacy between $(X,F)$ and $(X',F')$.
But for every $g\in G$
\begin{displaymath}
\phi_g\circ\tau^{-1}=(\tau\circ\phi_{g^{-1}})^{-1}
=(\sigma^G_{g^{-1}}\circ\tau)^{-1}
=\tau^{-1}\circ\sigma^G_g\;,
\end{displaymath}
thus
\begin{displaymath}
\sigma^G_g\circ F'=\tau\circ\phi_g\circ F\circ\tau^{-1}
=\tau\circ F\circ\phi_g\circ\tau^{-1}
=F'\circ\sigma^G_g\;;
\end{displaymath}
hence, $F'$ commutes with $\sigma^G$.
\iffalse
Let $c\in X'$ and let $x\in X$ satisfy $c=\tau(x)$: then
\begin{displaymath}
(F'(c))^g=(c_{F(x)})^g=c_{(F(x))^g}=c_{F(x^g)}
\end{displaymath}
and
\begin{displaymath}
F'(c^g)=(\tau\circ F)(\tau^{-1}((c_x)^g))
=(\tau\circ F)(\tau^{-1}(c_{x^g}))=\tau(F(x^g))
\end{displaymath}
with the rightmost terms of the identities being equal,
so $F'$ commutes with the natural right action of $G$ over $X'$.
\fi
By Hedlund's theorem, there exist a finite $\Neigh'\subseteq G$
and a map $f':A^{\Neigh'}\to A$ such that
$(F'(c))_g=f'\left(c^g|_{\Neigh'}\right)$
for all $c\in X'$, $g\in G$:
then $\left<X',\Neigh',f'\right>$
is a presentation of $(X,F)$ as a cellular automaton.
\qed
\end{proof}
The meaning of Theorem~\ref{thm_char}
is that $(X,F)$ has a CA presentation
with alphabet $A$ and tessellation group $G$,
if and only if $G$ can act on $X$ as it would do on $A^G$,
and without interfering with $F$.
This explains why the characterization
works for $CA(A,G)$, and not for $FCA(A,G)$:
the natural action, by itself,
is uncapable of telling the full shift from any other shift;
hence, any action on $X$ that ``emulates'' the natural action
shall not be able to tell whether $(X,F)$
has a presentation as a full CA or not.

Theorem~\ref{thm_char} has two immediate consequences.
The first one is a generalization, to our class of general CA,
of a principle first discovered by Hedlund~\cite{hed69} in dimension 1,
then extended by Richardson~\cite{ri72}
to classical CA of arbitrary dimension.
\begin{corollary} \label{cor_ca-inv}
Let $(X,F)\in CA(A,G)$.
If $F$ is bijective then $(X, F^{-1})\in CA(A,G)$.
\end{corollary}
\begin{proof}{}
Let $\phi$ be as in Theorem~\ref{thm_char}.
Then $X$ is discernible on $A$ by $\phi$, and
\begin{displaymath}
F^{-1}\circ\phi_g
=(\phi_{g^-1}\circ F)^{-1}
=(F\circ\phi_{g^{-1}})^{-1}
=\phi_g\circ F^{-1}
\end{displaymath}
for all $g\in G$.
Apply Theorem~\ref{thm_char}.
\qed
\end{proof}
Corollary~\ref{cor_ca-inv} can---and, in fact,
has been (cf.~\cite{ff03})---proved by purely topological means.
Our proof, however, gives some more hint
on the role of the tessellation group.
\begin{corollary}[Hedlund-Richardson's principle] \label{cor_HR}
Every bijective CA is reversible.
\end{corollary}
The second consequence of Theorem~\ref{thm_char}
is that existence of a presentation as CA actually depends
on the \emph{minimum number} of elements of the alphabet
and the \emph{isomorphism class} of the tessellation group:
which is intuitively true,
because isomorphic groups have ``isomorphic'' actions on equal spaces,
and because, if one has enough ``letters''
to be able to tell elements from each other via the action,
then having even more letters cannot be a bane.
\begin{proposition} \label{prop_A-G}
Let $A$ and $B$ be alphabets, and let $G$ and $\Gamma$ be f.g. groups.
\begin{enumerate}
\item \label{prop_A-G_p_AB}
If $|A|\leq|B|$ then $CA(A,G)\subseteq CA(B,G)$.
\item \label{prop_A-G_p_GGamma}
If $G$ is isomorphic to $\Gamma$ then $CA(A,G)=CA(A,\Gamma)$.
\end{enumerate}
\end{proposition}
\begin{proof}{}
To prove point~\ref{prop_A-G_p_AB},
let $\iota:A\to B$ be injective.
Let $(X,F)\in CA(A,G)$, and let $\phi$
satisfy point~\ref{thm_char_p_action} of Theorem~\ref{thm_char},
$\pi$ being the discerning map.
Then $X$ is discernible over $B$ by $\phi$,
$\iota\circ\pi$ being the discerning map.

To prove point~\ref{prop_A-G_p_GGamma},
let $\psi:G\to\Gamma$ be a group isomorphism.
Let $(X,F)\in CA(A,G)$ and let $\phi$
satisfy point~\ref{thm_char_p_action} of Theorem~\ref{thm_char},
$\pi$ being the discerning map.
Define $\phi'=\{\phi'_{\gamma}\}_{\gamma\in\Gamma}$ as
\begin{displaymath}
\phi'_{\gamma}=\phi_{\psi^{-1}(\gamma)}\;.
\end{displaymath}
It is straightforward to check that $\phi'$ is an action
which commutes with $F$.
Let $x_1\neq x_2$: if $g\in G$ is such that
$\pi(\phi_g(x_1))\neq\pi(\phi_g(x_2))$,
then $\pi(\phi'_{\psi(g)}(x_1))\neq\pi(\phi'_{\psi(g)}(x_2))$ as well.
Thus $\phi'$ satisfies condition 2 of Theorem~\ref{thm_char},
and $(X,F)\in CA(A,\Gamma)$.
From the arbitrariness of $(X,F)$ follows $CA(A,G)\subseteq CA(A,\Gamma)$:
by swapping the roles of $G$ and $\Gamma$
and repeating the argument with $\psi^{-1}$ in place of $\psi$
we obtain the reverse inclusion.
\qed
\end{proof}
%% end of Characterization of CA Dynamics via Group Actions

\section{Induced Subshifts} \label{sec_iss}

Let $X\subseteq A^G$ be a shift subspace.
We know that $X=\Xset^{A,G}_\Forb$ for some set $\Forb$ of patterns,
that is, $X$ is completely \textit{described} by $\Forb$
in the \textit{context} provided by $A$ and $G$.

Let now $\Gamma$ be a group having $G$ as a subgroup.
We want to define a new subshift $X'$ of $A^\Gamma$,
which is ``induced'' by $X$,
in the sense that $X'$ can be completely described by $X$.
But we had observed that $X$, in turn,
can be completely described by $\Forb$,
\emph{provided we know to be dealing with a subshift of $A^G$};
the first idea that comes to our mind,
is that $X'$ should then be completely described by $\Forb$ as well,
\emph{provided we know to be dealing with a subshift of $A^\Gamma$}.

This is precisely the content of
\begin{definition} \label{def_ind-s}
Let $X=\Xset^{A,G}_\Forb$ be a subshift,
and let $G\leq\Gamma$.
The \emph{subshift induced by $X$ on $A^\Gamma$} is
$X'=\Xset^{A,\Gamma}_\Forb$.
\end{definition}
\begin{example} \label{ex_gms}
Consider $A=\{0,1\}$, $G=\Zset$, $\Gamma=\Zset^2$, $\Forb=\{11\}$,
where $11$ is the pattern $p:\{0,1\}\subseteq\Zset\to A$
such that $p(0)=p(1)=1$.
Then $X=\Xset^{A,G}_\Forb$ is the \emph{golden mean shift}
(cf.~\cite{lm95});
a configuration $c:\Zset\to\{0,1\}$ belongs to $X$
if and only if it does not contain two adjacent 1's.
On the other hand, a configuration $\chi:\Zset^2\to\{0,1\}$
belongs to $X'=\Xset^{A,\Gamma}_\Forb$
if and only if no point on the square grid containing a 1
has his immediate right neighbour containing a 1 as well.
\qed
\end{example}
According to Definition~\ref{def_ind-s},
$X'$ is what we obtain instead of $X$,
when we interpret $\Forb$ as a description of a subshift
of $A^\Gamma$ instead of $A^G$,
that is, in the context provided by $\Gamma$ instead of $G$.

At first glance, Definition~\ref{def_ind-s} seems to be
a good solution to our ``subshift induction problem''.
However, we know from basic theory (cf.~\cite{lm95})
that different sets of patterns can define identical subshifts;
and we want induction to depend on
\emph{the object and not the description}.
We must then ensure that Definition~\ref{def_ind-s} is well posed
and $X'$ only depends on $X$ rather than $\Forb$,
\ie, $\Xset^{A,G}_{\Forb_1}=\Xset^{A,G}_{\Forb_2}$
must imply $\Xset^{A,\Gamma}_{\Forb_1}=\Xset^{A,\Gamma}_{\Forb_2}$.

In fact, we are going to discover much more.
We had noticed in Section~\ref{sec_defs}
that the image of a subshift via a UL-definable function
is a subshift;
thus, the most general form for a subshift of $A^G$ is
\begin{equation} \label{eq_ss-general}
X=F^{A,G}_f\left(\Xset^{A,G}_\Forb\right)\:,
\end{equation}
with $f:A^\Neigh\to A$, $\Neigh$ finite subset of $G$,
and $\Forb$ set of patterns with supports contained in G.
However, as we can choose to consider $\Forb$ as a description
of a subshift of either $A^G$ or $A^\Gamma$,
so we can choose to to consider $f$ as a description
of a UL-definable function on either $A^G$ or $A^\Gamma$.
It comes out that a more general thing we can check,
is the preservation of \emph{mutual inclusion}---instead
of just equality---between objects of the form (\ref{eq_ss-general}).

And this is precisely the content of
\begin{lemma} \label{lem_subg-UL}
Let $A$ be an alphabet,
and let $G$ and $\Gamma$ be f.g. groups with $G\leq\Gamma$.
For $i=1,2$, let $\Forb_i$ be a set of patterns on $A$
with supports contained in $G$,
let $\Neigh_i$ be a finite nonempty subset of $G$,
and let $f_i:A^{\Neigh_i}\to A$.
Then
\begin{displaymath}
F^{A,G}_{f_1}\left(\Xset^{A,G}_{\Forb_1}\right)
\subseteq F^{A,G}_{f_2}\left(\Xset^{A,G}_{\Forb_2}\right)
\;\;\mathrm{iff}\;\;
F^{A,\Gamma}_{f_1}\left(\Xset^{A,\Gamma}_{\Forb_1}\right)
\subseteq F^{A,\Gamma}_{f_2}\left(\Xset^{A,\Gamma}_{\Forb_2}\right)\;.
\end{displaymath}
\end{lemma}
\begin{proof}{}
Let $J$ be a set of representatives of the left cosets of $G$ in $\Gamma$
such that $1_G=1_\Gamma\in J$.
To simplify notation, we will write
\begin{displaymath}
X_i=\Xset^{A,G}_{\Forb_i}\;,\;
\Xi_i=\Xset^{A,\Gamma}_{\Forb_i}\;,\;
F_i=F^{A,G}_{f_i}\;,\;
\Phi_i=F^{A,\Gamma}_{f_i}\;,
\end{displaymath}
so that the thesis becomes
\begin{displaymath}
F_1(X_1)\subseteq F_2(X_2)\;\;\mathrm{iff}\;\;
\Phi_1(\Xi_1)\subseteq\Phi_2(\Xi_2)\;.
\end{displaymath}
For the ``if'' part, let $c\in F_1(X_1)$,
and let $x_1\in X_1$ satisfy $F_1(x_1)=c$.
Define $\xi_1\in A^\Gamma$ by $\xi_1(jg)=x_1(g)$
for all $j\in J$, $g\in G$:
then for all $j\in J$, $g\in G$, $p\in\Forb_1$
\begin{displaymath}
\restrict{\xi_1^{jg}}{\supp p}
=\restrict{x_1^g}{\supp p}
\neq p\;,
\end{displaymath}
hence $\xi_1\in\Xi_1$.
Put $\chi=\Phi_1(\xi_1)$:
by hypothesis, there exists $\xi_2\in\Xi_2$
such that $\Phi_2(\xi_2)=\chi$,
and by construction,
\begin{displaymath}
\chi(g)
=f_1(\restrict{\xi_1^g}{\Neigh_1})
=f_1((\restrict{x_1^g}{\Neigh_1})
=c(g)
\;\;\forall g\in G\;.
\end{displaymath}
Let $x_2=\restrict{\xi_2}{G}$:
then $x_2\in X_2$ by construction.
But
\begin{displaymath}
f_2(\restrict{x_2^g}{\Neigh_2})
=f_2(\restrict{\xi_2^g}{\Neigh_2})
=\chi(g)=c(g)
\;\;\forall g\in G\;,
\end{displaymath}
thus $c\in F_2(X_2)$.

For the ``only if'' part,
let $\chi\in\Phi_1(\Xi_1)$,
and let $\xi_1\in\Xi_1$
satisfy $\Phi_1(\xi_1)=\chi$.
For each $j\in J$, define $x_{1,j}\in A^G$
as $x_{1,j}(g)=\xi_1(jg)$ for all $g\in G$.
It is straightforward to check that $x_{1,j}\in X_1$ for all $j\in J$:
let $c_j=F_1(x_{1,j})$.
By hypothesis, for all $j\in J$
there exists $x_{2,j}\in X_2$
such that $F_2(x_{2,j})=c_j$:
define $\xi_2\in A^\Gamma$
by $\xi_2(jg)=x_{2,j}(g)$ for all $j\in J$, $g\in G$.
It is straightforward to check that $\xi_2\in\Xi_2$;
but for all $j\in J$, $g\in G$
\begin{displaymath}
f_2(\restrict{\xi_2^{jg}}{\Neigh_2})
=f_2(\restrict{x_{2,j}^g}{\Neigh_2})
=c_j(g)
=f_1(\restrict{x_{1,j}^g}{\Neigh_1})
=f_1(\restrict{\xi_1^{jg}}{\Neigh_1})
=\chi(jg)\;,
\end{displaymath}
thus $\chi\in\Phi_2(\Xi_2)$.
\qed
\end{proof}
The reason why Lemma~\ref{lem_subg-UL} is true, is the following.
Each left coset of $G$ can be thought of
as a ``slice'' of $\Gamma$ ``shaped'' as $G$.
If each pattern's support is contained in $G$,
then the constraint of not having a pattern in $\Forb_i$
can be applied either \emph{slice by slice}
or on the whole $\Gamma$ at once,
with the same results;
similarly, the neighbours of $\gamma$ w.r.t. $\Neigh_i$
will all belong to the same slice as $\gamma$,
so the $\Phi_i$'s can be made to operate
either slice by slice or on the whole $\Gamma$ at once,
with the same results.
This, however, means that the yes/no information
about the mutual inclusion of the $\Phi_i(\Xi_i)$'s
is deducible from the $\Forb_i$'s and the $f_i$'s alone,
and cannot be different from that on the $F_i(X_i)$'s.

Observe that the proof of Lemma~\ref{lem_subg-UL}
does not depend on the choice of the set $J$
of representatives of the left cosets of $G$ in $\Gamma$.
\begin{corollary} \label{cor_subg-UL}
In the hypotheses of Lemma~\ref{lem_subg-UL},
\begin{enumerate}
\item \label{cor_subg-UL_p_XtoFX}
\begin{math}
\Xset^{A,G}_{\Forb_1}
\subseteq F^{A,G}_{f_2}(\Xset^{A,G}_{\Forb_2})
\end{math}
iff
\begin{math}
\Xset^{A,\Gamma}_{\Forb_1}
\subseteq F^{A,\Gamma}_{f_2}(\Xset^{A,\Gamma}_{\Forb_2}),
\end{math}
\item \label{cor_subg-UL_p_FXtoX}
\begin{math}
F^{A,G}_{f_1}(\Xset^{A,G}_{\Forb_1})
\subseteq\Xset^{A,G}_{\Forb_2}
\end{math}
iff
\begin{math}
F^{A,\Gamma}_{f_1}(\Xset^{A,\Gamma}_{\Forb_1})
\subseteq\Xset^{A,\Gamma}_{\Forb_2},
\end{math}
and
\item \label{cor_subg-UL_p_XtoX}
\begin{math}
\Xset^{A,G}_{\Forb_1}
\subseteq\Xset^{A,G}_{\Forb_2}
\end{math}
iff
\begin{math}
\Xset^{A,\Gamma}_{\Forb_1}
\subseteq\Xset^{A,\Gamma}_{\Forb_2}.
\end{math}
\end{enumerate}
\end{corollary}
\begin{proof}
Consider the neighbourhood index $\{1_G\}$
and the local evolution function $f(1_G\mapsto a)=a$.
Apply Lemma~\ref{lem_subg-UL}.
\qed
\end{proof}
\begin{corollary} \label{cor_sofic-up}
Let $A$ be an alphabet,
let $G$ and $\Gamma$ be f.g. groups with $G\leq\Gamma$,
and let $\Forb$ be a set of patterns on $A$
with supports contained in $G$.
If $\Xset^{A,G}_\Forb$ is sofic then $\Xset^{A,\Gamma}_\Forb$ is sofic.
\end{corollary}
\begin{proof}
By hypothesis, $\Xset^{A,G}_\Forb=F(\Xset^{A,G}_{\Forb'})$
for some UL-definable function $F$ and finite set of patterns $\Forb'$.
Apply points~\ref{cor_subg-UL_p_XtoFX}
and~\ref{cor_subg-UL_p_FXtoX}
of Corollary~\ref{cor_subg-UL}.
\qed
\end{proof}
Now, under the same hypotheses on $G$, $\Gamma$, $A$, and $\Forb$,
suppose $X'=\Xset^{A,\Gamma}_\Forb$ is sofic.
This means that
a finite set $\Forb'$ of patterns \emph{over $\Gamma$}
and a function $f':A^{\Neigh'}\to A$
with $\Neigh'$ finite subset \emph{of $\Gamma$},
such that
\begin{math}
X'=F^{A,\Gamma}_{f'}\left(\Xset^{A,\Gamma}_{\Forb'}\right).
\end{math}
Is it then possible for $X=\Xset^{A,G}_\Forb$ \emph{not} to be sofic?
In fact, the finitary description for $X'$ provided by $\Forb'$ and $f'$
takes advantage of the (at least, \textit{a priori}) greater complexity
of the group $\Gamma$ w.r.t. the group $G$;
however, it is also true that $\Forb$ alone
yields enough information to describe $X'$
in the context provided by $\Gamma$.
It would not be surprising, then,
if the information provided by $\Forb$
in the context provided by $G$
yielded enough information to describe $X$;
we state this as a conjecture,
hoping it will be either proved or disproved in the future.
\begin{conjecture} \label{conj-sofic-down}
With the notation of Corollary~\ref{cor_sofic-up},
suppose $\Xset^{A,\Gamma}_\Forb$ is sofic.
Then $\Xset^{A,G}_\Forb$ is sofic.
\end{conjecture}
%% end of Induced Subshifts

\section{Induced Cellular Automata} \label{sec_ica}
In the previous section,
we have learned to construct a subshift on a group
from a subshift on a subgroup;
while doing this, we have received some insight
on how the underlying mechanism can also work for UL-definable functions.
It then comes to our mind that similar mechanics
could be applied to another field where locality is the key factor,
that is, the field of cellular automata.
This time, we can give our definition
after having already done the bulk of the work.
\begin{definition} \label{def_induced}
Let $\Acal=\left<X,\Neigh,f\right>$ be a CA
with alphabet $A$ and tessellation group $G$,
and let $\Gamma$ be a f.g. group such that $G\leq\Gamma$.
The CA \emph{induced by $\Acal$ on $\Gamma$}
is the cellular automaton
\begin{equation} \label{eq_ca-ind}
\Acal'=\left<X',\Neigh,f\right>\;,
\end{equation}
where $X'$ is the subshift induced by $X$ on $A^\Gamma$.
\end{definition}
Observe how Lemma~\ref{lem_subg-UL}
ensures that $\Acal'$ is well defined.
\begin{example} \label{ex_ind-ca}
Let $A=\{0,1\}$, $G=\Zset$, $\Gamma=\Zset^2$, $\Neigh=\{-1,1\}$,
$f(-1\mapsto x,1\mapsto y)=x+y-2xy$:
then $\Acal=\left<A^\Zset,\Neigh,f\right>$
is \emph{Wolfram's rule 90},
such that the next value of each point
is the exclusive OR of the current values
of its leftmost and rightmost neighbours.
The same rule applies to $\Acal'=\left<A^{\Zset^2},\Neigh,f\right>$,
which can be seen as the joining
of infinitely many copies of $\Acal$
along a vertical line.
\qed
\end{example}
Definition~\ref{def_induced} is similar to the one
given in~\cite{csms99} for CA over the full shift;
ours, however, works for the broader class of CA over subshifts.
(We still have, however, the constraint on finite alphabets,
which~\cite{csc08} tries to overcome at least for the full shift.)
As in the case of the induced subshift---which, by the way,
is the support of the induced CA---$\Acal'$ is what we obtain
by interpreting the local descriptions given by $\Forb$, $\Neigh$, and $f$,
in the context provided by $\Gamma$ instead of $G$.

It must be remarked that, in general, $\Acal'$ is not conjugate to $\Acal$.
For instance, if $G$ is a proper nontrivial subgroup
of a finite group $\Gamma$,
then there can be no bijection between $A^G$ and $A^\Gamma$,
let alone conjugacies of cellular automata.

On the other hand, it had already been shown in~\cite{csms99}
that, in the case of CA over full shifts,
some important properties---notably, surjectivity---are preserved
in the passage from the original CA to the induced one;
which is not surprising, because intuitively
$F_f^{A,\Gamma}$ is going to operate ``slice by slice'' on $A^\Gamma$,
each ``slice'' being ``shaped'' as $G$.
The next statement extends the aforementioned result
from the case $X=A^G$ to the general case when $X$ is an arbitrary subshift.
\begin{theorem} \label{thm_induced-surj}
Let $\Acal=\left<X,\Neigh,f\right>$ be a CA
with alphabet $A$ and tessellation group $G$,
let $G\leq\Gamma$,
and let $\Acal'$ be the CA induced by $\Acal$ on $\Gamma$.
\begin{enumerate}
\item \label{thm_induced-surj_p_surj}
$\Acal$ is surjective iff $\Acal'$ is surjective.
\item \label{thm_induced-surj_p_preinj}
$\Acal$ is preinjective iff $\Acal'$ is preinjective.
\item \label{thm_induced-surj_p_inj}
$\Acal$ is injective iff $\Acal'$ is injective.
\end{enumerate}
\end{theorem}
\begin{proof}
Let $\Forb$ satisfy $X=\Xset^{A,G}_\Forb$ (and $X'=\Xset^{A,\Gamma}_\Forb$).
Take $J$ as in proof of Lemma~\ref{lem_subg-UL}.

To prove the ``if'' part of point~\ref{thm_induced-surj_p_surj},
suppose $\Acal$ has a GoE pattern $p$.
By contradiction, assume that there exists $\chi\in X'$
such that $\restrict{F_{\Acal'}(\chi)}{\supp p}=p$.
Let $c$ be the restriction of $\chi$ to $G$.
Then, since both $\Neigh$ and $\supp p$ are subsets of $G$ by hypothesis,
\begin{displaymath}
(F_\Acal(c))(x)=f\left(\restrict{c^x}{\Neigh}\right)
=f\left(\restrict{\chi^x}{\Neigh}\right)
=(F_{\Acal'}(\chi))(x)=p(x)
\end{displaymath}
for every $x\in\supp p$:
this is a contradiction.

To prove the ``only if'' part of point~\ref{thm_induced-surj_p_surj},
suppose $\Acal'$ has a GoE pattern $\pi$.
By hypothesis, there exists $\chi\in X'$
such that $\restrict{\chi}{\supp\pi}=\pi$.
For all $j\in J$ define $c_j\in A^G$ as
\begin{displaymath}
c_j(g)=\chi(jg)\;\;\forall g\in G\;,
\end{displaymath}
and for all $j\in J$ such that $jG\cap\supp\pi\neq\emptyset$
define the pattern $p_j$ over $G$ as
\begin{displaymath}
p_j(x)=\pi(jx)\,\;\forall x\;\mathrm{s.t.}\;jx\in\supp\pi\;.
\end{displaymath}
Observe that $c_j\in X$ for all $j$,
and that $p_j=\restrict{c_j}{jG\cap\supp\pi}$ when defined.
But at least one of the patterns $p_j$ must be a GoE for $\Acal$:
otherwise, for all $j\in J$, either $jG\cap\supp\pi=\emptyset$,
or there would exist $k_j\in X'$
such that $\restrict{F_\Acal(k_j)}{\supp p_j}=p_j$.
In this case, however, $\kappa\in A^\Gamma$ defined by
$\kappa(jg)=k_j(g)$ for all $j\in J$, $g\in G$
would satisfy $\kappa\in X'$
and $\restrict{F_{\Acal'}(\kappa)}{\supp\pi}=\pi$,
against $\pi$ being a GoE for $\Acal'$.

For the ``if'' part of point~\ref{thm_induced-surj_p_preinj},
suppose $c_1,c_2\in X$ differ on all and only
the points of a finite nonempty $U\subseteq G$,
but $F_\Acal(c_1)=F_\Acal(c_2)$.
For all $j\in J$, $g\in G$,
put $\chi_1(jg)=c_1(g)$,
and set $\chi_2(jg)$ as $c_2(g)$ if $j=1_\Gamma$, $c_1(g)$ otherwise.
Then $\chi_1$ and $\chi_2$ belong to $X'$
and differ precisely on $U$.
Moreover, for every $\gamma\in\Gamma$,
either $\gamma\in G$ or $\gamma\Neigh\cap G=\emptyset$,
so either
\begin{math}
(F_{\Acal'}(\chi_i))(\gamma)=(F_\Acal(c_i))(\gamma)
\end{math}
or
\begin{math}
(F_{\Acal'}(\chi_1))(\gamma)=(F_{\Acal'}(\chi_2))(\gamma).
\end{math}

For the ``only if'' part of point~\ref{thm_induced-surj_p_preinj},
suppose $\Acal$ is preinjective.
Let $\chi_1,\chi_2\in X'$ differ on all and only
the points of a finite nonempty $U'\subseteq\Gamma$.
For $i\in\{1,2\}$, $\gamma\in\Gamma$, let $c_{i,\gamma}$
be the restriction of $\chi_i^\gamma$ to $G$:
these are all in $X$,
because a pattern occurring in $c_{i,\gamma}$ also occurs in $\chi_i$,
and cannot belong to $\Forb$.
Let
\begin{math}
U_\gamma=\{g\in G\mid c_{1,\gamma}(g)\neq c_{2,\gamma}(g)\}:
\end{math}
then $|U_\gamma|\leq|U|$ for all $\gamma\in\Gamma$,
plus $U_\gamma\neq\emptyset$ for at least one $\gamma$.
For such $\gamma$, there exists $g\in G$ such that
\begin{math}
(F_\Acal(c_{1,\gamma}))(g)\neq(F_\Acal(c_{2,\gamma}))(g):
\end{math}
then by construction
\begin{math}
(F_{\Acal'}(\chi_1))(\gamma g)
\neq(F_{\Acal'}(\chi_2))(\gamma g)
\end{math}
as well.

The proof of point~\ref{thm_induced-surj_p_inj}
is straightforward to see.
For the ``if'' part, let $c_1\neq c_2$, $F_\Acal(c_1)=F_\Acal(c_2)$,
and consider $\chi_i(\gamma)=c_i(g)$ iff $\gamma=jg$.
For the ``only if'' part, given $\chi_1\neq\chi_2$,
consider $c_{i,j}(g)=\chi_i(jg)$,
and observe that $F_\Acal(c_{1,j})\neq F_\Acal(c_{2,j})$
for at least one $j\in J$.
\qed
\end{proof}
The reason why Theorem~\ref{thm_induced-surj} is true,
is similar to the one given for Lemma~\ref{lem_subg-UL}:
the global evolution function of the induced CA
operates ``slice by slice'' on the support of the induced CA;
this, however, is the induced subshift,
and is already ``sliced'' suitably for $F_{\Acal'}$.
Moreover, each of the listed global properties
can be expressed in local terms:
for instance, surjectivity is equivalent to absence of GoE patterns,
even in our broader context (cf.~\cite{ff03}).
Pay attention, however, that these properties
are usually r.e. or co-r.e., but not computable.

Observe that, as in the proof of Lemma~\ref{lem_subg-UL},
the choice of $J$ is arbitrary.
\begin{example} \label{ex_ind-properties}
Let $\Acal$ be as in Example~\ref{ex_ind-ca}.
It is a good exercise in cellular automata theory
to check that each configuration has exactly four predecessors
according to $\Acal$,
that is, for every $c:\Zset\to A$
there exist four distinct $c_i:\Zset\to A$
such that $F_\Acal(c_i)=c$.
(Hint: fix four patterns $00,01,10,11$.)
Thus $\Acal$ is surjective, but not injective;
Theorem~\ref{thm_induced-surj} then says that $\Acal'$
is also surjective and noninjective.
\qed
\end{example}
Surjectivity and preinjectivity are always shared by $\Acal$ and $\Acal'$,
even when these two properties are not equivalent:
This fact was used in~\cite{csms99}
to prove that Moore-Myhill's theorem does not hold
for full CA with tessellation group
containing a free subgroup on two generators,
starting from two counterexamples on the free group $\Free_2$.

Having observed that $\Acal$ and $\Acal'$ may well be non-conjugate,
we are left with a different question:
is it possible to \emph{embed} the original CA into the induced one?
After all, we have \emph{kept} the same local descriptions,
and \emph{enlarged} the group,
so we should expect the induced dynamics
to be \emph{richer} than the original.
Moreover, since the global evolution function of $\Acal'$
operates slice by slice,
we should expect that,
after having \emph{fixed a point on each slice},
we should be able to reproduce $\Acal$ into $\Acal'$.

And this is precisely the content of
\begin{lemma} \label{lem_ind-emb}
Let $A$ be an alphabet,
and let $G$ and $\Gamma$ be f.g. groups with $G\leq\Gamma$;
let $\Acal=\left<X,\Neigh,f\right>$
be a CA with alphabet $A$ and tessellation group $G$,
and let $\Acal'=\left<X',\Neigh,f\right>$
be the CA induced by $\Acal$ over $\Gamma$.
Let $J$ be a set of representatives of the left cosets of $G$ in $\Gamma$,
and let $\iota_J:A^G\to A^\Gamma$ be defined by
\begin{equation} \label{eq_iotaJ}
(\iota_J(c))(\gamma)=c(g)\;\;\mathrm{iff}\;\;\exists j\in J:\gamma=jg\;.
\end{equation}
Then $\iota_J$ is an embedding of $\Acal$ into $\Acal'$, so that
\begin{equation} \label{eq_iotaJ-ca}
\iota_J(\Acal)=\left<\iota_J(X),\Neigh,f\right>
\end{equation}
is a CA conjugate to $\Acal$.
In particular, $CA(A,G)\subseteq CA(A,\Gamma)$.
\end{lemma}
\begin{proof}{}
First, we observe that $\iota_J$ is injective and $\iota_J(X)\subseteq X'$.
In fact, if $c_1(g)\neq c_2(g)$,
then $(\iota_J(c_1))(jg)\neq(\iota_J(c_2))(jg)$ for all $j\in J$.
Moreover, should a pattern $p$ exist
such that $(\iota_J(c))(\gamma x)=p(x)$ for all $x\in\supp p\subseteq G$,
by writing $\gamma=jg$ and applying (\ref{eq_iotaJ})
we would find $c(gx)=p(x)$ for all $x\in\supp p$,
a contradiction.

Next, we show that $\iota_J$ is continuous.
Let $S$ be a f.s.o.g. for $G$,
$\Sigma$ a f.s.o.g. for $\Gamma$.
Let $R\geq 0$, and let
\begin{displaymath}
E_R=\{g\in G\mid \exists j\in J\mid jg\in D^\Gamma_{R,\Sigma}\}\;.
\end{displaymath}
Since the writings $\gamma=jg$ are unique and $D^\Gamma_{R,\Sigma}$ is finite,
$E_R$ is finite too.
Let $E_R\subseteq D^G_{r,S}$:
if
\begin{math}
\restrict{c_1}{D^G_{r,S}}=\restrict{c_2}{D^G_{r,S}},
\end{math}
then
\begin{math}
\restrict{\iota_J(c_1)}{D^\Gamma_{R,\Sigma}}
=\restrict{\iota_J(c_2)}{D^\Gamma_{R,\Sigma}}.
\end{math}

Next, we show that $\iota_J$ is a morphism of d.s.
For every $c\in A^G$, $\gamma=jg\in\Gamma$, $x\in\Neigh$
we have $\gamma x\in jG$ and $(\iota_J(c))(\gamma x)=(\iota_J(c))(jgx)=c(gx)$.
Thus,
\begin{displaymath}
((F_{\Acal'}\circ\iota_J)(c))(\gamma)
=f\left(\iota_J(c)^\gamma|_\Neigh\right)
=f\left(c^g|_\Neigh\right)
=(F_\Acal(c))(g)=((\iota_J\circ F_\Acal)(c))(\gamma)\;,
\end{displaymath}
so that $F_{\Acal'}\circ\iota_J=\iota_J\circ F_\Acal$.
Moreover,
\begin{math}
F_{\Acal'}(\iota_J(X))=\iota_J(F_\Acal(X))\subseteq\iota_J(X)
\end{math}
because $F_\Acal(X)\subseteq X$.

Finally, we observe that $\iota_J(X)$ is a subshift.
In fact, if $X=\Xset^{A,G}_\Forb$,
then $\iota_J(X)=\Xset^{A,\Gamma}_{\Forb\cup\Forb'}$, where
\begin{equation} \label{eq_moreforb}
\Forb'=\left\{p\in A^{\{j_1g,j_2g\}}\mid
j_1,j_2\in J, g\in G, j_1\neq j_2,
p(j_1g)\neq p(j_2g)\right\}\;.
\end{equation}
It is straightforward that
$\iota_J(X)\subseteq\Xset^{A,\Gamma}_{\Forb\cup\Forb'}$.
Let $\chi\in\Xset^{A,\Gamma}_{\Forb\cup\Forb'}$:
then $c(g)=\chi(jg)$ is well defined, and $\chi=\iota_J(c)$ by construction.
Moreover, for every $g\in G$, $p\in F$, and any $j\in G$
\begin{math}
(c^g)_{\supp p}=(\chi^{jg})_{\supp p}\neq p,
\end{math}
so $c\in X$ and $\chi\in\iota_J(X)$.
\qed
\end{proof}
Observe that, in the hypotheses of Lemma~\ref{lem_ind-emb},
$\iota_J$ depends explicitly on $J$,
so it may, in general, show ``better'' or ``worse'' properties
according to the choice of $J$:
such properties, however, have no effect
on the \emph{abstract dynamics} of $\iota_J(\Acal)$,
which is always the same as $\Acal$'s.

Moreover, we are \emph{not} assuming $1_\Gamma\in J$;
hence, in general, $E_R\not\subseteq D^G_{R,S}$,
even if $S\subseteq\Sigma$.
\begin{example} \label{ex_counter}
Let $\Gamma=\Zset^2$,
$G=\{(x,0)\mid x\in\Zset\}$,
$S=\{(1,0)\}$, $\Sigma=\{(1,0),(0,1)\}$, and
\begin{displaymath}
J=\{(1,0)\}\cup\{(0,y)\mid y\in\Zset,y\neq 0\}\:.
\end{displaymath}
Then $E_1=\{(0,0),(-1,0),(-2,0)\}\not\subseteq D^G_{1,S}$.
\qed
\end{example}
Lemma~\ref{lem_ind-emb} says that growing the tessellation group
does not shrink the class of presentable dynamics.
This fact and Proposition~\ref{prop_A-G} together yield
\begin{theorem} \label{thm_incl}
Let $A$, $B$ be alphabets and $G$, $\Gamma$ be f.g. groups.
If $|A|\leq|B|$ and $G$ is isomorphic to a subgroup of $\Gamma$,
then $CA(A,G)\subseteq CA(B,\Gamma)$.
\end{theorem}
\begin{proof}{}
Let $G\cong H\leq\Gamma$. Then
\begin{math}
CA(A,G)=CA(A,H)\subseteq CA(A,\Gamma)\subseteq CA(B,\Gamma).
\end{math}
\qed
\end{proof}
\begin{corollary} \label{cor_Free}
Let $\Free_n$ be the free group on $n<\infty$ generators.
For every alphabet $A$ and every $n>1$, $CA(A,\Free_n)=CA(A,\Free_2)$.
\end{corollary}
\begin{proof}{}
Clearly, every $\Free_n$ with $n>1$
has a free subgroup on two generators:
because of Theorem~\ref{thm_incl},
$CA(A,\Free_2)\subseteq CA(A,\Free_n)$.
However, it is a well-known fact in group theory
(cf.~\cite{mks76}, Section 2.4, Problem 2)
that $\Free_2$ has a free subgroup on infinitely many generators,
hence also has a free subgroup on $n$ generators for every $n>0$:
because of Theorem~\ref{thm_incl},
$CA(A,\Free_n)\subseteq CA(A,\Free_2)$.
\qed
\end{proof}
Observe that the inclusion in one direction
also works for $n=1$, with $\Free_1=\Zset$.
Since Moore-Myhill's theorem does not hold for the latter class
(cf.~\cite{csms99})
we know that $FCA(\{0,1\},\Zset)\neq FCA(\{0,1\},\Free_2)$;
however, a similar statement for the corresponding $CA$-classes
has not been established yet.
In fact, the structure of $\Free_2$
is intrinsically much more complex than that of $\Zset$,
where the same cannot be said of the other $\Free_n$'s,
which somehow ``contain each other''.
\begin{conjecture} \label{conj_ca-free-z}
$CA(A,\Zset)\neq CA(A,\Free_2)$.
\end{conjecture}
Now, if we look at (\ref{eq_moreforb}),
the set of ``additional constraints'' $\Forb'$
seems a bit cumbersome.
Why is it necessary to take note of all the pairs $(j_1g,j_2g)$?
Can we only make the checks on the pairs $(j_1,j_2)$,
and use the smaller set
\begin{displaymath}
\Forb''=\left\{p\in A^{\{j_1,j_2\}}\mid
j_1,j_2\in J,p(j_1)\neq p(j_2)\right\}\;,
\end{displaymath}
or can we?

The problem with the idea of replacing $\Forb'$ with $\Forb''$,
is that we are forgetting that $\Gamma$ can be noncommutative:
thus, $j\gamma$ is not forced to equal $\gamma j$,
which is what we get when we try to check
whether the configuration $\chi$ has a pattern with support $\{j_1,j_2\}$.
On the other hand, if $\Forb''$ were always a good choice,
then, for $\iota_J(X)$ to be of finite type,
it would suffice $X$ being of finite type
and $G$ being of finite index in $\Gamma$,
independently on the choice of $J$.
This seems just too good to be true;
and, in fact, it is false.
\begin{theorem} \label{thm_ft-non-emb}
Let $\Gamma$ be the group of ordered pairs $(i,k)$,
$i\in\{0,1\}$, $k\in\Zset$ with the product
\begin{displaymath}
(i_1,k_1)(i_2,k_2)=(i_1+i_2-2i_1i_2,(-1)^{i_2}k_1+k_2)\;.
\end{displaymath}
Let $A=\{a,b\}$, $G=\{(0,k),k\in\Zset\}\leq\Gamma$, and $J=\{(0,0),(1,0)\}$.
Then $\iota_J(A^G)$ is not a shift of finite type.
\end{theorem}
\begin{proof}{}
Let $S=\{(1,0),(0,1)\}$:
it is straightforward to check that $\left<S\right>=\Gamma$.

By contradiction, assume
that $\iota_J(A^G)=\Xset^{A,\Gamma}_\Forb$ with $|\Forb|<\infty$;
it is not restrictive to choose $\Forb$ so that
$\supp p=D^\Gamma_{M,S}$ for all $p\in\Forb$.
Let $\delta\in A^\Gamma$ satisfy $\delta(x)=b$ iff $x=(0,0)$:
then $\delta\not\in\iota_J(A^G)$,
so there must exist $p\in\Forb$, $\eta\in\Gamma$
such that $\delta^{\eta}|_{\supp p}=p$.
It is straightforward to check that there exists
exactly one $y\in D^\Gamma_M$ such that $p(y)=b$,
and that $y=\eta^{-1}=(i,(-1)^{1-i}x)$ if $\eta=(i,x)$.

Now, for all $k\in\Zset$ we have
$d^\Gamma_S((0,k),(1,k))=\|(1,2k)\|^\Gamma_S=2|k|+1$.
This can be checked by observing the following two facts.
Firstly, $(1,2k)=(1,0)(0,t)\ldots(0,t)$, with $2|k|$ factors $(0,t)$,
and $t=1$ or $t=-1$ according to $k>0$ or $k<0$.
Secondly, multiplying $(i,x)$ on the right by $(0,1)$ or $(0,-1)$
does not change the value of $i$,
while multiplying $(i,x)$ on the right by $(1,0)$
does not change $|x|$:
hence, at least one multiplication by $(1,0)$
and $2|k|$ multiplications by either $(0,1)$ or $(0,-1)$
are necessary to reach $(1,2k)$ from $(0,0)$.

For $i\in\{0,1\}$ let $\gamma_i=(i,2M+1)$.
Let $\chi\in A^\Gamma$ be such that $\chi(\gamma)=b$
iff $\gamma=\gamma_0$ or $\gamma=\gamma_1$:
then $\chi\in\iota_J(A^G)$.
However, since $\eta^{-1}\in D^\Gamma_{M,S}$,
for all $x\in D^\Gamma_M(\eta^{-1})$
we have $\gamma_0\eta x\in D^\Gamma_{2M}(\gamma_0)$.
Hence, either $x=\eta^{-1}$, $\gamma_0\eta x=\gamma_0$,
and $\chi^{\gamma_0\eta}(x)=b$;
or $x\neq\eta^{-1}$,
$0<d_S(\gamma_0,\gamma_0\eta x)\leq 2M<4M+3=d_S(\gamma_0,\gamma_1)$,
and $\chi^{\gamma_0\eta}(x)=a$.
Thus, $\restrict{\chi^{\gamma_0\eta}}{\supp p}=p$:
this is a contradiction.
\qed
\end{proof}
The reason why Theorem~\ref{thm_ft-non-emb} is true,
is that, in general,
one cannot get an upper bound on $d_S(j_1g,j_2g)$
only by looking at $d_S(j_1,j_2)$:
this happens because the product is made \emph{on the wrong side}.
Consequently, one should not expect to determine
finitely many constraints on the $jg$'s,
only from finitely many constraints on the $j$'s.
\begin{corollary} \label{cor_ft-nondyn}
For cellular automata on arbitrary f.g. groups,
finiteness of type is not invariant by conjugacy.
In particular, for subshifts on arbitrary f.g. groups,
finiteness of type is not a topological property.
\end{corollary}
The first statement in Corollary~\ref{cor_ft-nondyn}
seems to collide with Theorem 2.1.10 of~\cite{lm95},
stating that any two conjugate subshifts of $A^\Zset$
are either both of finite type or both not of finite type.
Actually, in the cited result, conjugacies are always intended
as being between \emph{shift dynamical systems},
which is a much more specialized situation than ours.
Moreover, the tessellation group is always $\Zset$,
so that \emph{the action is also the same},
while we have different groups and different actions.
Last but not least, translations are UL-definable
if and only if the translating factor is \emph{central},
\ie, commutes with every element in the tessellation group:
thus, the only groups where \emph{all} the translations are UL-definable
are the abelian groups.
On the other hand, the second statement remarks the well-known phenomenon
that homeomorphisms do not preserve finiteness of type,
not even in ``classical'' symbolic dynamics.
\begin{example} \label{ex_evsh}
Let $\Forb=\{10^{2n+1}1\mid n\in\Nset\}$:
the subshift $X=\Xset^{\{0,1\},\Zset}_\Forb$
is called the \emph{even shift}.
It can be proved (cf.~\cite{lm95}, Section 3.1)
that $X$ is not a shift of finite type.
However, $X$ is
\begin{enumerate}
\item non-empty---it contains the configuration with all 0's,
\item compact---as a subshift,
\item metrizable---with the distance inherited by the full shift,
\item totally disconnected---because the full shift is, and
\item perfect---every point of $X$ can be seen as the limit point
of some sequence of elements of $X$
taking the value 1 only finitely many times.
\end{enumerate}
By a theorem of Brouwer,
the even shift is homeomorphic to the Cantor set,
thus also to the full shift---which \emph{is} of finite type.
\qed
\end{example}
%%For instance, the \emph{even shift}
%%(there is always an even number of 0's between any two 1's)
%%is not of finite type,
%%but it is homeomorphic to the Cantor set, thus also to the full shift.

Our attempt at finding a general criterion for finiteness of type
having crashed against an unsurmountable obstacle,
we choose to switch our aim towards a more modest target.
What if we add conditions on the \emph{way} $G$ is related to $\Gamma$,
and are more careful in the choice of $J$?

A possible answer is given by
\begin{theorem} \label{thm_dp-ft-ok}
Let $H$ and $K$ be f.g. groups;
let $S$ be a finite set of generators for $H$
such that $1_H\not\in S$ and $H=\left<S\right>$;
let $\Gamma=H\times K$, $G=\{1_H\}\times K$, $J=H\times\{1_K\}$.
Let $A$ be an alphabet and let
\begin{displaymath}
\Forb_S=\left\{p\in A^{\{(1_H,1_K),(s,1_K)\}}\mid
s\in S\cup S^{-1}\setminus\{1_H\},
p((1_H,1_K))\neq p((s,1_K))\right\}\;.
\end{displaymath}
For every set $\Forb$ of patterns on $A$
with supports contained in $G$,
\begin{math}
\iota_J(\Xset^{A,G}_\Forb)=\Xset^{A,\Gamma}_{\Forb\cup\Forb_S}.
\end{math}
In particular, if $X\subseteq A^G$ is a shift of finite type,
then $\iota_J(X)$ is also a shift of finite type.
\end{theorem}
\begin{proof}{}
First, observe that $\Forb_S\subseteq\Forb'$,
where $\Forb'$ is given by (\ref{eq_moreforb}), so that
$\iota_J(\Xset^{A,G}_\Forb)=\Xset^{A,\Gamma}_{\Forb\cup\Forb'}
\subseteq\Xset^{A,\Gamma}_{\Forb\cup\Forb_S}$.

Let now $\chi\in A^\Gamma\setminus\iota_J(X)$;
suppose that no $p\in\Forb$ occurs in $\chi$.
Let $h_1,h_2\in H$, $k\in K$ satisfy $\chi((h_1,k))\neq\chi((h_2,k))$,
and let $h_1^{-1}h_2=s_1s_2\cdots s_N$
be a writing of minimal length of the form (\ref{eq_gen}).
For $i\in\{0,\ldots,N\}$ let
\begin{math}
a_i=\chi(h_1s_1\ldots s_i,k)\;;
\end{math}
for $i\in\{1,\ldots, N\}$ define
\begin{math}
p_i:\{(1_H,1_K),(s_i,1_K)\}\to A
\end{math}
by
\begin{math}
p_i(1_H,1_K)=a_{i-1}
\end{math}
and
\begin{math}
p_i(s_i,1_K)=a_i\;.
\end{math}
Since $a_0\neq a_N$, $a_{i-1}\neq a_i$ for some $i$:
then $p_i\in\Forb_S$ and
\begin{math}
\restrict{\chi^{(h_1s_1\cdots s_{i-1},k)}}{\supp p_i}=p_i.
\end{math}
Since $\chi$ is arbitrary,
\begin{math}
\Xset^{A,\Gamma}_{\Forb\cup\Forb_S}\subseteq\iota_J(\Xset^{A,G}_\Forb).
\end{math}
\qed
\end{proof}
The reason why Theorem~\ref{thm_dp-ft-ok} holds
is that, even if $\Forb_S$ puts less restraints than $\Forb'$,
it is also true that the components $J$ and $G$ of the direct product
\emph{do not interfere with each other} in the multiplication,
so any $j_1g_1\cdots j_ng_n$
can be rewritten as $j_1\cdots j_ng_1\cdots g_n$ without any trouble,
and with the result still being of the form $jg$, $j\in J$, $g\in G$.

Observe that, for Theorem~\ref{thm_dp-ft-ok} to hold,
$G$ needs not to be of finite index in $\Gamma$.
However, the other hypotheses are quite strong,
especially the ones on the structure of $\Gamma$
as a direct product with $G$ as a factor:
this seems unlikely to be improved easily
because of Theorem~\ref{thm_ft-non-emb},
where $\Gamma$ is a \emph{semi-direct product} with $G$ as a factor.
\begin{example} \label{ex_dp-ft-ok}
Let $A=\{0,1\}$;
let $H$ and $K$ be two distinct copies of $\Zset$
with $S=\{1\}$.
Identify $\Gamma=H\times K$ with $\Zset^2$,
$J=H\times\{0\}$ and $G=\{0\}\times K$ with $\Zset$.
Let $p:\{(0,0),(0,1)\}\to\{0,1\}$
satisfy $p(0,0)=p(0,1)=1$:
then $X=\Xset^{A,G}_{\{p\}}$ can be identified with the golden mean shift.
Let $p_{01},p_{10}:\{(1,0),(1,1)\}\to\{0,1\}$ be defined by
\begin{displaymath}
p_{01}(1,0)=0\;,\;p_{01}(1,1)=1\;,\;p_{10}(1,0)=1\;,\;p_{10}(1,1)=0\;.
\end{displaymath}
Then
\begin{math}
\iota_J(X)=\Xset^{A,\Gamma}_{\left\{p,p_{01},p_{10}\right\}}.
\end{math}
\qed
\end{example}
We conclude with a statement about sofic shifts.
\begin{theorem} \label{thm_ft-sofic}
Let $A$, $G$, $\Gamma$, and $J$ be as in Lemma~\ref{lem_ind-emb}.
Suppose $\iota_J(X)$ is a shift of finite type
for every shift of finite type $X\subseteq A^G$.
Then $\iota_J(X)$ is a sofic shift for every sofic shift $X\subseteq A^G$.
\end{theorem}
\begin{proof}{}
Let $X=F(Y)$ for some shift of finite type $Y\subseteq A^G$
and UL-definable function $F:A^G\to A^G$.
Let $\Neigh\subseteq G$, $|\Neigh|<\infty$, and $f:A^\Neigh\to A$
be such that $(F(c))_g=f(c^g|_\Neigh)$ for all $c\in A^G$, $g\in G$;
let $\Acal=\left<A^G,\Neigh,f\right>$
and let $F'$ be the global evolution function of $\iota_J(\Acal)$.
By Lemma~\ref{lem_ind-emb}, $F'\circ\iota_J=\iota_J\circ F$,
so that $\iota_J(X)=\iota_J(F(Y))=F'(\iota_J(Y))$
is the image of a shift of finite type via a UL-definable function.
\qed
\end{proof}
We suspect that the hypotheses in Theorem~\ref{thm_ft-sofic}
are, in fact, redundant.
Again, we state this as a conjecture,
hoping for it to be proved or disproved in the future.
\begin{conjecture} \label{conj_sofic}
Let $A$, $G$, $\Gamma$, and $J$ be as in Lemma~\ref{lem_ind-emb}.
Suppose $X\leq A^G$ is a sofic shift.
Then $\iota_J(X)$ is a sofic shift.
\end{conjecture}
%% end of Induced Cellular Automata

\section{Conclusions} \label{sec_concl}
At the end of our trek, we have seen
how to get CA presentations of dynamical systems,
and to construct new shift subspaces and cellular automata
by enlarging their underlying groups.
We have also remarked the properties of old objects
inherited by the new ones,
and taken note of some exceptions.
Finally, we have observed how enlarging the group
makes the class of presentable dynamics grow.

However, there is surely much work to do.
In particular, much to our shame,
we were not able to either prove or disprove
Conjectures~\ref{conj-sofic-down} and~\ref{conj_sofic},
nor to determine whether they have found a solution.
Also additional conditions on the discerning action $\phi$
in the proof of Theorem~\ref{thm_char},
such to get a characterization of full CA dynamics,
has been painfully missed.

Aside of looking ourselves for the answers to such questions,
our hope is that our modest work can be of interest,
or even use, to researchers in the field.
%% end of Conclusions

\section{Acknowledgements} \label{sec_ack}
The author was partially supported by the project
``The Equational Logic of Parallel Processes''
(nr.~060013021) of The Icelandic Research Fund.
We also thank Luca Aceto, Anna Ing\'olfsd\'ottir,
Tommaso Toffoli, Patrizia Mentrasti,
Tullio Ceccherini-Silberstein,
Carlos Mart\'{\i}n-Vide, Jarkko Kari,
and the the anonymous referees for the LATA 2008 conference
for their many suggestions and encouragements.
%% end of Acknowledgements


\begin{thebibliography}{99}

\bibitem{c04}
S. Capobianco.
\textit{Structure and Invertibility in Cellular Automata}.
PhD thesis, University of Rome ``La Sapienza'', 2004.

\bibitem{c08lata}
S. Capobianco.
Induced Subshifts and Cellular Automata.
To appear in issue of \textit{Lect.Not. Comp. Sci.}
dedicated to LATA 2008 conference proceedings.

\bibitem{csc08}
T.G. Ceccherini-Silberstein, M. Coornaert.
Induction and restriction of cellular automata.
To appear on \textit{Erg. Th. Dyn. Syst.}
Preprint: \verb+www-irma.u-strasbg.fr/~coornaer/induction.pdf+

\bibitem{csms99}
T. G. Ceccherini-Silberstein, A. Mach\`\i, F. Scarabotti.
Amenable groups and cellular automata.
\textit{Ann. Inst. Fourier, Grenoble} \textbf{42} (1999), 673--685.

\bibitem{ff03}
F. Fiorenzi.
Cellular automata and strongly irreducible shifts of finite type.
\textit{Theor. Comp. Sci.} \textbf{299} (2003), 477--493.

\bibitem{hed69}
G. A. Hedlund.
Endomorphisms and automorphisms of the shift dynamical system.
\textit{Math. Syst. Th.} \textbf{3} (1969), 320--375.

\bibitem{lm95}
D. Lind, B. Marcus.
\textit{An introduction to symbolic dynamics and coding.}
Cambridge University Press 1995.

\bibitem{mm93}
A. Mach\'{\i}, F. Mignosi.
Garden of Eden Configurations
for Cellular Automata on Cayley Graphs on Groups.
\textit{SIAM J. Disc. Math.} \textbf{6} (1993), 44--56.

\bibitem{mks76}
W. Magnus, A. Karrass, D. Solitar.
\textit{Combinatorial Group Theory. Presentations of Groups in Terms of Generators and Relations.}
Dover Publications, Inc., 1976.

\bibitem{mo62}
E. F. Moore.
Machine models of self-reproduction.
\textit{Proc. Symp. Appl. Math.} \textbf{14} (1962), 17--33.

\bibitem{my62}
J. Myhill.
The converse of Moore's Garden-of-Eden theorem.
\textit{Proc. Am. Math. Soc.} \textbf{14} (1962), 685--686.

\bibitem{ri72}
D. Richardson.
Tessellations with local transformations.
\textit{J. Comp. Syst. Sci.} \textbf{6} (1972), 373--388.

\end{thebibliography}
\end{document}